\documentclass[preprint,12pt]{elsarticle}

\usepackage{xspace}

 \usepackage[latin1]{inputenc}
 \usepackage[T1]{fontenc}
 \usepackage{graphicx}
 \usepackage{caption}
 \usepackage{subcaption}
 \usepackage{url}
 \usepackage{alltt}
 \usepackage[english]{babel}
 \usepackage{amsmath}
 \usepackage{amssymb}
 \usepackage{tikz}
 \usepackage{pgfplotstable}
 \usepackage{pgfplots}
  \usepackage{float}

 \newcommand{\Insitu}{In Situ\xspace}
 \newcommand{\intransit}{in transit\xspace}
 
 \newcommand{\codesaturne}{{\fontfamily{ppl}\fontshape{it}\selectfont Code\_Saturne}\xspace}
\newcommand{\melissa}{Melissa\xspace}
 
 \newcommand{\zmq}{\textsc{ZeroMQ}\xspace}

\newcommand{\bigO}				{\textit{O}}
\newcommand{\hide}[1]{}

\journal{Simulation Modelling Practice and Theory}
\journal{xxxx}

\begin{document}

\begin{frontmatter}



\title{Large scale \intransit  computation of quantiles for  ensemble runs}


\author{Alejandro Rib\'es$^1$, Th\'eophile Terraz$^2$, Bertrand Iooss$^{3,4}$, Yvan Fournier$^3$, Bruno Raffin$^2$}

\address{$^1$ EDF Lab Paris-Saclay, France\\
$^2$ Univ. Grenoble Alpes, Inria, CNRS, Grenoble INP, LIG, France\\
$^3$ EDF Lab Paris-Chatou, France\\
$^4$ Institut de Math\'ematiques de Toulouse, Universit\'e Paul Sabatier, Toulouse, France
}

\begin{abstract}


The classical approach for quantiles computation requires availability of the full sample before ranking it. 
In uncertainty quantification of numerical simulation models, this approach is not suitable at exascale as large ensembles of simulation runs would need to gather a prohibitively large amount of data.  
This problem is solved thanks to an on-the-fly and iterative approach based on the Robbins-Monro algorithm.
This approach relies on Melissa, a file avoiding, adaptive, fault-tolerant and elastic framework. 
On a validation case producing $11$ TB of data, which consists in $3000$ fluid dynamics parallel simulations on a $6$M cell mesh, it allows on-line computation of spatio-temporal maps of percentiles.

\end{abstract}


\end{frontmatter}



\section{Introduction}
\label{sec:introduction}

On-line analytics for large scale numerical simulations has demonstrated its potential to contain the I/O bottleneck, improve simulation and analytics performance, and overall reduce the human time scientists spent in handling large data sets. 
A key enabler is the availability of one-pass analytics algorithms, i.e. algorithms that  can operate on a reduced  window of data recently produced  by the simulation to avoid the need for massive storage and its associated performance pitfalls.
In this paper, we focus on ensemble runs where multiple instances, usually thousands, of the simulation are run to sweep across the parameter space. 
Such process is typically led in engineering practice during the uncertainty quantification stage of some simulation models~\cite{smi14,morbal14}. 
This stage mainly consists in computing statistical quantities of the model outputs (mean, variance, quantile, probability of threshold exceedance, \ldots) or estimating sensitivity indices between model outputs and inputs (linear correlation coefficients, Sobol' indices, \ldots)~\cite{baudin2016open,ansmar14}.

A major difficulty arises when the ensemble runs produce massive amount of data that have to be statistically aggregated, making them extremely vulnerable to the I/O bottleneck.  
To keep a manageable amount of data,  the  classical approach, used in most of the studies, consists in reducing the resolution of the simulation and the sampling points for the statistics~\cite{ghahig17}. 
A more suitable technique would be to use one-pass statistical algorithms.
One-pass statistical algorithms, also  called iterative, online or even parallel statistics have the interesting property of requiring only to store the current  results that can next be updated  with incoming new samples.      
One-pass variance algorithms were proposed in~\cite{welford1962note,chan1982updating,Finch-incremental-stats-2009}. 
Numerically   stable,  one-pass   formulas   for  arbitrary   centered statistical     moments      and     co-moments      are     presented in~\cite{bennett2009numerically,pebay2008}.
\cite{pebay2008}   also contains update formulas for  higher order moments (skewness, kurtosis and  more), setting  the  bases  for  a module  of  parallel statistics      in      the     VTK      scientific      visualization toolkit~\cite{pebay2011design}.
In   this   context   the   one-pass algorithms enables  to compute partial  results in parallel  before to perform  a  reduction  to  get   the  final  result.  
These  iterative statistics were used for computing large scale parallel statistics for a     single    simulation     run     either     from    raw     data files~\cite{bennett2011featureTVCG},          compressed          data files~\cite{lakshminarasimhan2011isabela}             or            in situ~\cite{bennett2012combining}.     
In sensitivity analysis of model outputs, for the estimation of Sobol' indices,  \cite{gilarn17}  introduced a one-pass iterative computation for the case of a scalar model output, while \cite{terraz-melissa-SC17} applied the iterative covariance formulas on massive output data (a spatio-temporal model output).

Various   packages  are designed for managing uncertainty quantification and sensitivity analysis  from ensemble runs (see for example \cite{ghahig17,baudin2016open}).
However, they all rely on classical non-iterative algorithms, requiring to accumulate first all simulation results in file or memory if doable. 
Based on a different  architecture, the Melissa framework~\cite{isav-16,terraz-melissa-SC17} has been recently proposed for  the  on-line data aggregation of high resolution  ensemble runs. 
Other in situ processing  frameworks (see \cite{glean08,Damaris:2012,FlexIO13,ccgrid:14,Marrinan-streaming-JCS-2019}) enable in situ and in transit processing but for a  single simulation.
In Melissa, each simulation handles  its output as soon as available   to a set of staging nodes. 
These nodes  process these incoming data   to update the statistics on a first-come first-served basis thanks to the one-pass algorithm.  
This in transit processing mode enables to fully avoid  storage of intermediate data on disks. 
Melissa runs  a parallel server that  stores the current state  of the computed statistics.  
Simulations dynamically connects to  this server and send their  results a soon as available to  update the statistics. 
This  architecture,  complemented  with a  fault  tolerance  mechanism allows for  efficient elastic executions.  

Melissa  currently supports the estimation of the following statistical quantities: standard deviation, skewness, kurtosis, minimum, maximum, threshold exceedance probability and Sobol' indices.   
This paper focuses on extending the Melissa framework to the quantile estimation issue, useful quantity for risk analysis, outlier detection or computation of non-parametric confidence intervals.
Low or high-order quantiles are often required during uncertainty quantification studies, especially for industrial safety issues \cite{derdev08,mor12,morbal16,ioomar18}.
Standard approaches deal with the problem of quantile estimation of scalar outputs \cite{gly96,hesnel98,cangar08}.
However, simulation models often return more complex objects as outputs, such as temporal curves and spatial fields, which can be considered as functions.
Recent studies have considered quantiles of one-dimensional functional outputs (temporal curves) \cite{popioo13,ribpou14,nanhel16} and demonstrated the interest for the practitioners to compute these functional quantiles.

Quantiles being order-statistics, the straightforward and classical approach consists in ordering the sample, then finding the appropriate quantiles~\cite{davnag03}. This strategy necessitates the storage of the full-sample, thus making its application in an on-line  context impossible. 
To avoid this storage, one can use stochastic algorithms which are devoted to the recursive estimation of statistical quantities.
For instance, the Robbins-Monro algorithm allows for the iterative computation of quantiles~\cite{robmon51,duf97}.
It has been introduced in the context of simulation model by \cite{kohkry14}, but only for a scalar output. 
In this paper, a new robust version of this algorithm is developed.
Then, instead of providing quantiles for a limited sample of probes as usually done, full spatio-temporal model outputs can be considered in order to estimate ubiquitous (i.e. everywhere in space and time) multidimensional and time varying quantiles.

The paper is organized as follows. After presenting the quantile estimation issue and proposed algorithm (Sec.~\ref{sec:quantiles}),  its implementation in the   \melissa architecture  is discussed (Sec.~\ref{sec:melissa}) before presenting the   experimental results (Sec.~\ref{sec:experiments}).  The visual analysis of the obtained data reveals the potential of ubiquitous quantile statistics  (Sec.~\ref{sec:visu}).   A conclusion closes the paper   (Sec.~\ref{sec:conclusion}).

\section{In Transit Computation of Quantiles }\label{sec:quantiles}

Quantiles are important in descriptive statistics because they display variation in samples of a statistical population without making any assumptions of the underlying statistical distribution. In the context of ensemble analysis and visualization, a priori knowledge about output data distributions is in general not known thus making these statistics highly useful. In fact, quantiles have long been  used for data visualization and understanding, the earliest example being the Tukey boxplot~\cite{tukey1977exploratory}, which is a method for graphically depicting scalar data distribution through their quartiles. In the present work, we are interested in the non-parametric characterization of the output variability of ensemble runs. We remark that these outputs not being scalar values, but ubiquitous multidimensional and time varying quantiles, their  computation, visualization and interpretation represents a challenge.



%
\hide{\theophile{
We consider a numerical simulation on a fixed mesh. The simulation computes values for different fields $u$ for each mesh
element $x$.
The simulation progresses in time through various time steps $t$ (all simulations simulate the same time steps).
The same simulation runs $N$ times with different input parameters.
Let call $i$ the $i^{th}$ simulation.
Our goal is to compute statistics over all simulations for each mesh element at each time step and for each field  $u(i,x,t)$.
}}


Computing statistics from $N$ samples classically requires  $\bigO(N)$ memory space to store  these samples.  
But if the statistics can be {\it computed in one-pass} (also called
iterative, online   or even parallel~\cite{pebay2008}), i.e. if the
current value can be updated as soon as a new sample is available, the
memory requirement goes down to  $\bigO(1)$ space.   With this approach, not only simulation results do not need to be
saved, but they can be consumed in any order, loosening synchronization constraints on the simulation executions.

Estimation procedures based on recursive algorithms are efficient techniques that are able to deal with large and voluminous samples.
The Robbins-Monro algorithm \cite{robmon51} is such a procedure and has been developed in many situations and applications \cite{benmet90,duf97,kusyin03}.
When the variable under study is not a simple scalar but a functional variable, the literature is much less abundant.
One can cite the works \cite{smayao06,carcen13} where infinite dimensional Banach or Hilbert space are considered, but such developments remain preliminary.
Our present work considers a high-dimensional vector of scalar variables (coming from the discretized spatio-temporal field) that are treated independently of each other.
Dealing with the functional space where the spatio-temporal field lives remains a challenge and will be the topic of further works.




\subsection{Empirical Quantile Estimator}

Let us consider a $N$-sample $(Y_1,\ldots,Y_N)$ of i.i.d. random variables from an unknown distribution $f_Y(y)$. 
We look for an estimator $\hat{q}_{\alpha}$ of the  $\alpha$-quantile $q_\alpha$ defined by:
\begin{equation}
\mathbb{P}( Y \leq q_\alpha) = \alpha \;,
\end{equation}
which is sometimes written as
\begin{equation}
q_\alpha = \inf\{y|\mathbb{P}( Y \leq y) \geq \alpha\} \;.
\end{equation}

The classical estimator of the $\alpha$-quantile is the empirical quantile, based on the notion of order statistics \cite{davnag03}.
Essentially, we associate with the sample $(Y_1,\ldots,Y_N)$ the ordered sample $(Y_{(1)}, \ldots, Y_{(N)})$ in which  $Y_{(1)} \leq \ldots \leq Y_{(N)}$. 
The empirical estimator then writes
\begin{equation}
\label{quantileempirique}
\hat{q}_{\alpha} = Y_{(\lfloor \alpha N \rfloor +1)} ,
\end{equation}
where $\lfloor x \rfloor$ is the integer part of $x$.
When the probability density of $Y$ is differentiable in $y_\alpha$, a central limit theorem ($N\to \infty$) exists, shown for example in \cite{davnag03}, which says that $\hat{q}_{\alpha}$  is an asymptotically normal estimator with variance $\alpha(1-\alpha)/[(N+2) f_Y^2(y_\alpha)]$.


\subsection{Tuning a Robbins-Monro Estimator Algorithm}

In the particular case of a quantile estimation, the Robbins-Monro estimator \cite{robmon51} consists in updating the quantile estimate  at each new observation $Y_n$ with the following rule
\newcommand{\1}{1 \hspace*{-0.4ex}\rule{0.10ex}{1.5ex}\hspace{0.4ex}} 
\begin{equation}\label{eq:RM}
q_\alpha(n+1) = q_\alpha(n) - \frac{C}{n^\gamma}\left( \1_{Y_{n+1 \leq q_\alpha(n)}} - \alpha \right) \;,
\end{equation}
with $n =1\ldots N$, $q_\alpha(n)$ the $\alpha$-quantile estimate at the $n$th observation, $q_\alpha(1) = Y_1$ an independent realization of $Y$, $\1_{x}$ the indicator function, $C$ a strictly positive constant and $\gamma \in ]0,1]$ the step of the gradient descent of the stochastic algorithm.
When $\gamma \in ]0.5,1]$ and under several hypotheses, this algorithm has been shown to be consistent and asymptotically normal.
The Robbins-Monro algorithm has been introduced in the context of simulation model (scalar) output by \cite{kohkry14} and has been used by \cite{labgam16} to solve the problem of conditional quantile estimation of stochastic simulation models.

In the equation (\ref{eq:RM}), $C$ and $\gamma$ have to be chosen.
In the following, we fix $C=1$ and concentrate our efforts to tune the $\gamma$ values.
Indeed, in the applications of this paper, the choice $C=1$ will be satisfactory because the variables under study have variations of the order of unity (the quantile updates involved in Eq. (\ref{eq:RM}) will have the good order of magnitude).
A finer tuning of this constant in practical applications will be developed in a future study.

Asymptotically (for $N\to \infty$), a value $\gamma=1$ is known to be optimal.
However, in practical studies, $N$ is not large enough to reach the asymptotic regime.
For the type of engineering studies we consider ($Y$ is the output of a costly computer code, see for example \cite{cangar08,ioomar18}), $\alpha=0.95$ and $N$ is in the order of several hundreds of simulated values.
To understand the algorithm behavior, a first numerical test is performed with $N=1000$, $\alpha=0.95$ and $Y$ following a standard Gaussian distribution ${\cal N}(0,1)$.
$N=1000$ is the typical order of magnitude of our studies (our application in Section \ref{sec:experiments} will use $3000$ simulations).

\begin{figure}[!ht]
\begin{subfigure}[t]{0.49\linewidth}
  \centering
  \captionsetup{width=.8\linewidth}
  \includegraphics[width=1\linewidth]{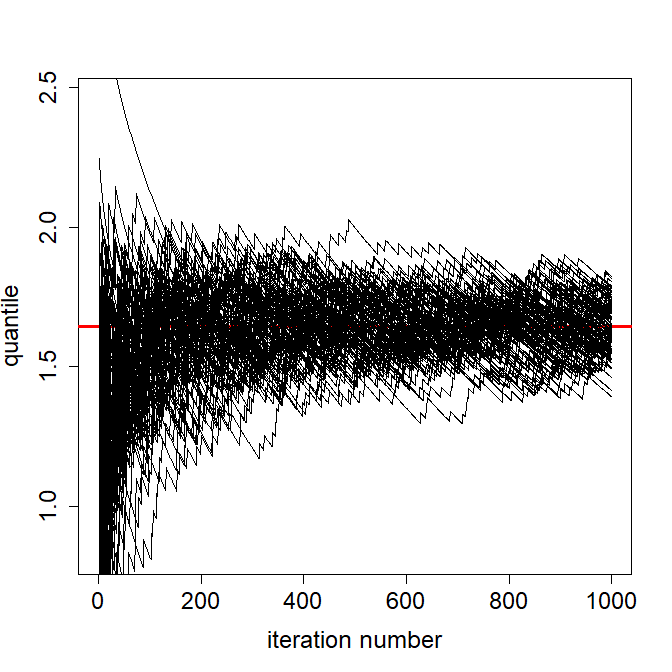}
  \caption{$\gamma=0.5$.}
  \label{subfig:gamma05}
\end{subfigure}
\begin{subfigure}[t]{0.49\linewidth}
  \captionsetup{width=.8\linewidth}
  \includegraphics[width=1\linewidth]{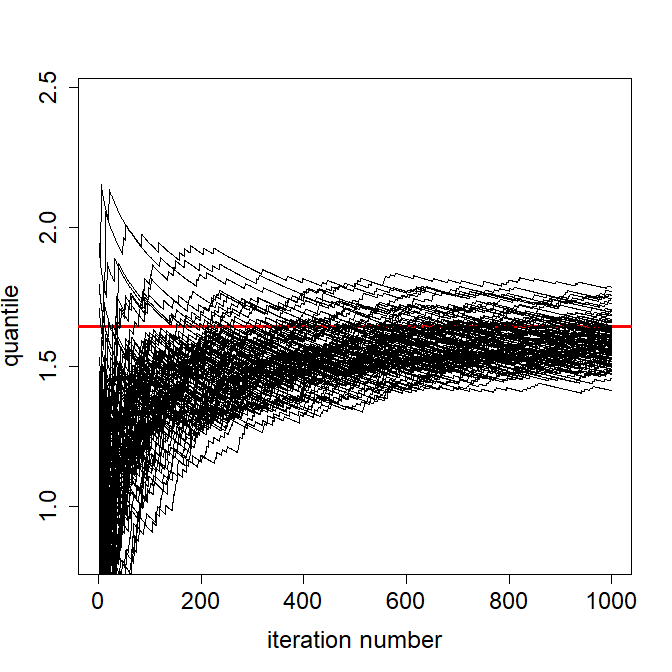}
  \caption{$\gamma=0.7$.}
  \label{subfig:gamma07}
\end{subfigure}

\begin{subfigure}[t]{0.49\linewidth}
  \centering
  \includegraphics[width=1\linewidth]{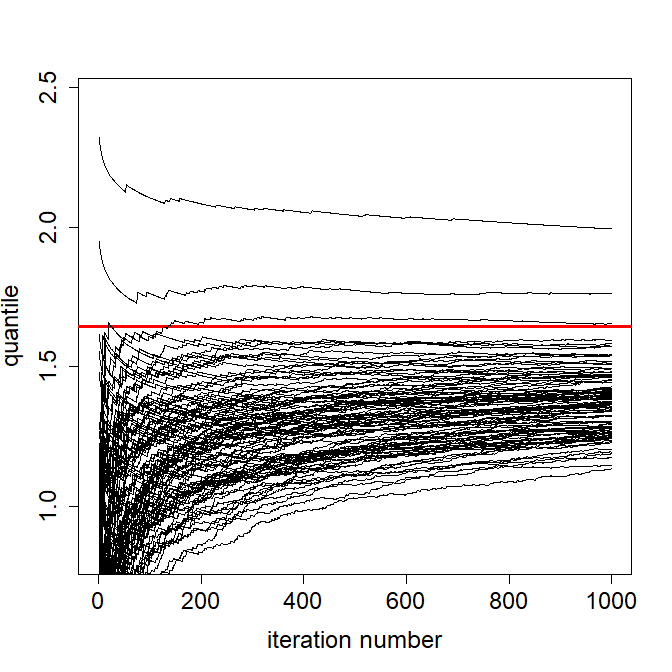}
  \caption{$\gamma=0.9$.}
  \label{subfig:gamma09}
\end{subfigure}
\begin{subfigure}[t]{0.49\linewidth}
  \includegraphics[width=1\linewidth]{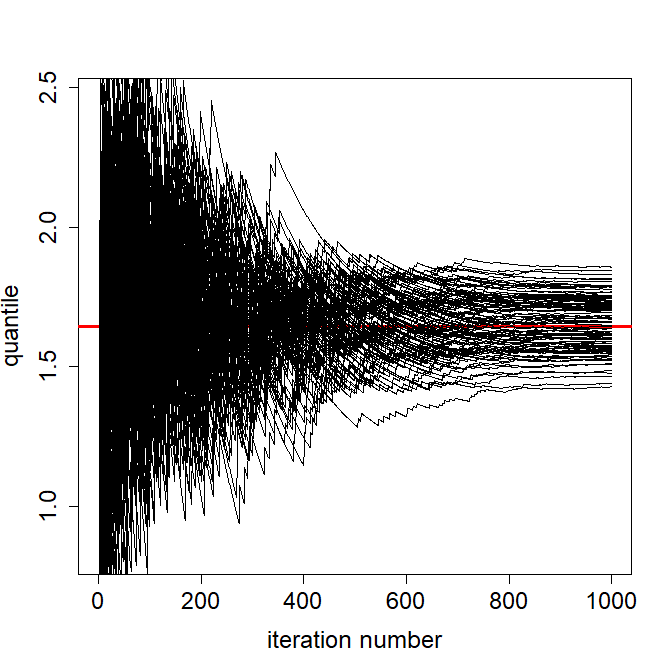}
  \caption{Linear $\gamma$ evolution.}
  \label{subfig:gammalinear}
\end{subfigure}

\vspace{0.1cm}
  \caption{For different choices of $\gamma$, simulation of $100$ independent trajectories ($n=1,\ldots,1000=N$) of the Robbins-Monroe estimation of the $0.95$-quantile of a ${\cal N}(0,1)$ variable. The red horizontal line is the exact $0.95$-quantile.} 
\vspace{-0.3cm}
  \label{fig:RMconvergence}
\end {figure}

The Figure \ref{fig:RMconvergence} shows $100$ different and independent trajectories of the Robbins-Monro quantile estimates $q_{0.95}(n)$ for $n=1,\ldots,1000$ considering different values of $\gamma$.
One can observe that:
\begin{itemize}
    \item small values of $\gamma$ (fig. \ref{subfig:gamma05}) induce large mixing of the quantile estimate during its evolution (all along the iterations). Then, the problem is a lack of stabilization because the final number of iterations could not be large enough;
    \item larger values of $\gamma$ (figs. \ref{subfig:gamma07} and \ref{subfig:gamma09}) induce small perturbations during the evolution of the quantile estimate. Then, the problem arrives when the initialization value ($n=1$) is far from the  quantile exact value because the evolution cannot correct it enough.
\end{itemize}
The Figure \ref{subfig:gammalinear} shows the result of an algorithm (that will be detailed below) which consists in having small values of $\gamma$ during the first iterations and large values of $\gamma$ during the last iterations.
A good convergence seems to have been reached (well-centered and not too dispersed distribution of the values around the exact value).

Another issue is our needs of robustness for the choice of $\gamma$.
Indeed, we look for $\gamma$ values which can work for different distributions of $Y$ (which are unknown in practice).
The problem is that good $\gamma$ values for a certain type of probability distribution give bad results for another type of distribution (for example, $\gamma=0.6$ gives good results for a normal distribution and incorrect quantile estimates for a uniform one). 
Therefore, we introduce a new way to deal with the Robbins-Monro algorithm by defining $\gamma$ as a function of $n$.
The heuristic formula, inspired by a linear temperature profile choice in the simulated annealing algorithm, is the following:
\begin{equation}
\gamma(n) = 0.1 + 0.9 \frac{n-1}{N-1} \;.
\end{equation}
The idea is to have a strong mixing properties at the beginning of the algorithm (with small $\gamma$), then to slow down the potential variation of the quantile estimation all along the iterations of the algorithm.

\begin{figure}[!ht]
\begin{subfigure}[t]{0.49\linewidth}
  \centering
  \captionsetup{width=.8\linewidth}
  \includegraphics[width=1\linewidth]{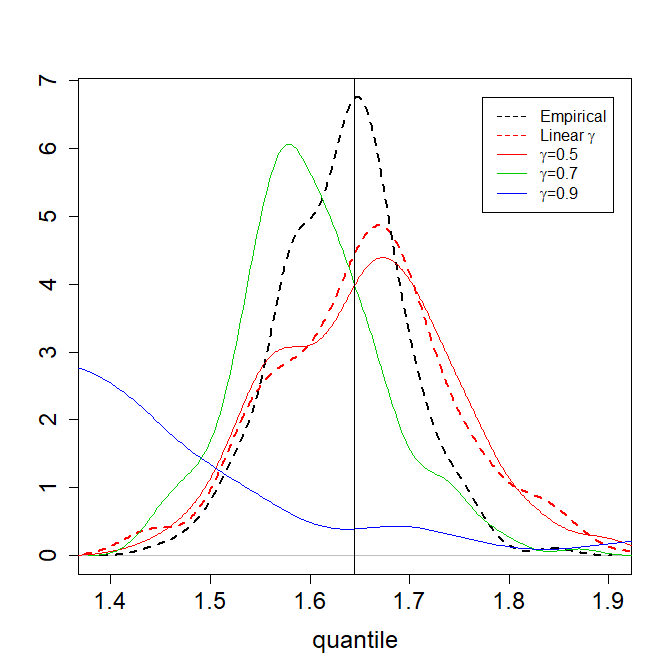}
  \caption{Gaussian, $Y\sim {\cal N}(0,1)$.}
  \label{subfig:RMtests_norm}
\end{subfigure}
\begin{subfigure}[t]{0.49\linewidth}
  \captionsetup{width=.8\linewidth}
  \includegraphics[width=1\linewidth]{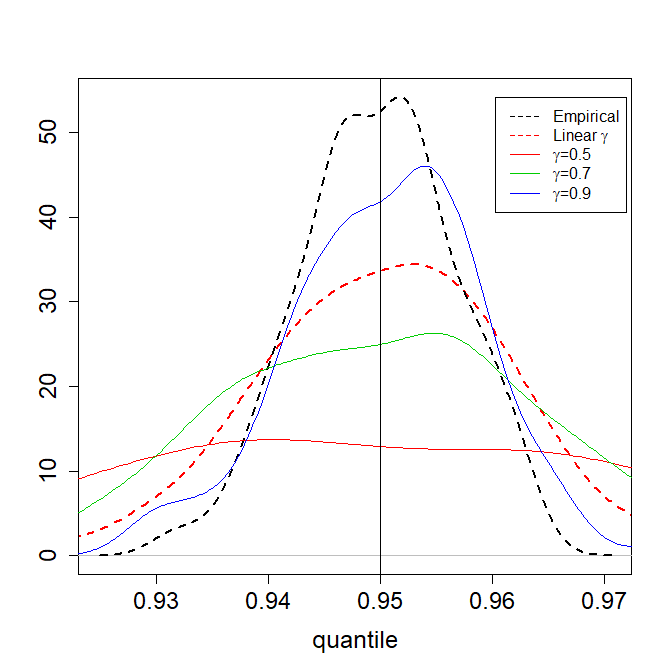}
  \caption{Uniform, $Y\sim {\cal U}[0,1]$.}
  \label{subfig:RMtests_unif}
\end{subfigure}

\begin{subfigure}[t]{0.49\linewidth}
  \centering
  \includegraphics[width=1\linewidth]{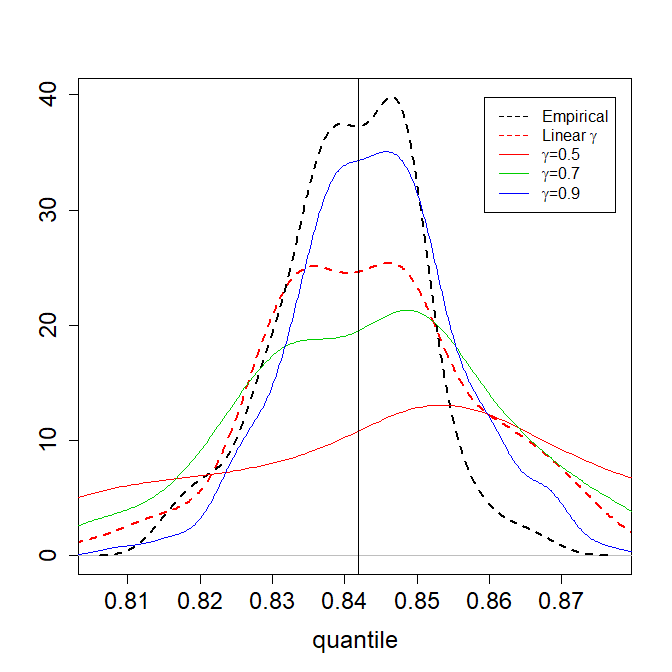}
  \caption{Triangular, $Y\sim {\cal T}(0,0.5,1)$.}
  \label{subfig:RMtests_tria}
\end{subfigure}
\begin{subfigure}[t]{0.49\linewidth}
  \includegraphics[width=1\linewidth]{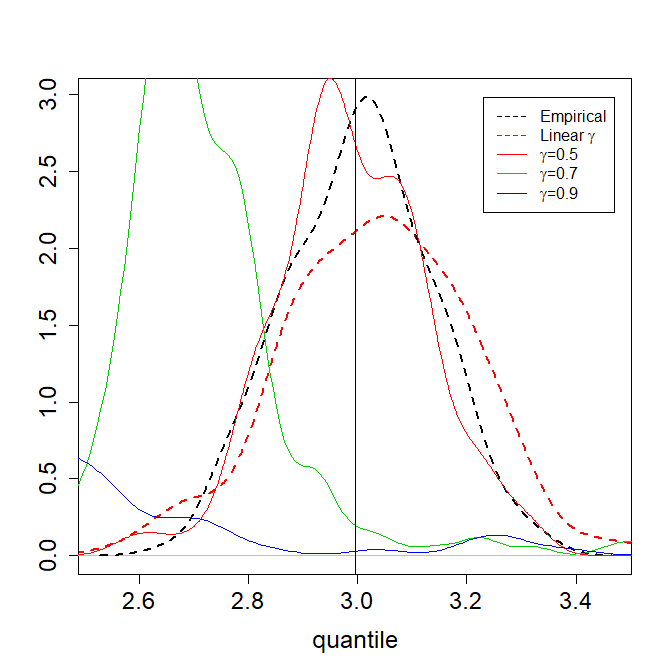}
  \caption{Exponential, $Y\sim {\cal E}(1)$.}
  \label{subfig:RMtests_exp}
\end{subfigure}

\vspace{0.1cm}
  \caption{For different choices of probability density functions of $Y$, probability densities of different estimators of the $0.95$-quantile of $Y$: Empirical estimator and Robbins-Monro estimators (with $\gamma=0.5$, $\gamma=0.7$, $\gamma=0.9$ and the linear profile). The vertical line is the exact $0.95$-quantile.} 
\vspace{-0.3cm}
  \label{fig:RMtests}
\end {figure}

Several numerical tests on different distributions of $Y$ and simple analytical functions (where the true quantile can be known) have been performed  to calibrate and validate this linear $\gamma$-profile.
Figure \ref{fig:RMtests} shows the results of four tests considering different probability density functions for $Y$, $\alpha=0.95$ and $N=1000$.
Comparisons are made between the results given by the linear $\gamma$-profile and by different constant values for $\gamma$.
We are particularly interested by knowing if we can obtain similar results with the Robbins-Monro estimator to those of the empirical estimator, which is our reference (because it is based on the storing of all sample values that is not the case with the Robbins-Monro estimators).
For each estimator, the test consists in repeating $200$ times the algorithm  to obtain distributed values of the estimates.
So, distribution-based comparisons are made.

For all the cases (figs. \ref{subfig:RMtests_norm}, \ref{subfig:RMtests_unif}, \ref{subfig:RMtests_tria} and \ref{subfig:RMtests_exp}), the distributions of estimates obtained by the linear $\gamma$-profile are  well-centered, not too dispersed and rather close to the empirical estimator based results.
This corresponds to the robustness we look for.
We observe also that a good constant $\gamma$ value in some cases is a  really poor choice in other cases.
For example, $\gamma=0.5$ gives excellent results in  figs. \ref{subfig:RMtests_norm} and \ref{subfig:RMtests_exp} and dramatic ones in figs. \ref{subfig:RMtests_unif} and \ref{subfig:RMtests_tria}, while  $\gamma=0.6$ gives excellent results in figs. \ref{subfig:RMtests_unif} and \ref{subfig:RMtests_tria} and poor ones in  figs. \ref{subfig:RMtests_norm} and \ref{subfig:RMtests_exp}.
All these results confirm the choice of the linear $\gamma$-profile that will be used in our practical study in the following.

\section{The Melissa Framework}\label{sec:melissa}

We present in this section an overview of the Melissa framework (see~\cite{terraz-melissa-SC17} for a more detailed description).

\subsection{\melissa Architecture}

\begin{figure}[!ht]
  \centering
  \vspace{-1.5cm}
  \includegraphics[width=1\linewidth]{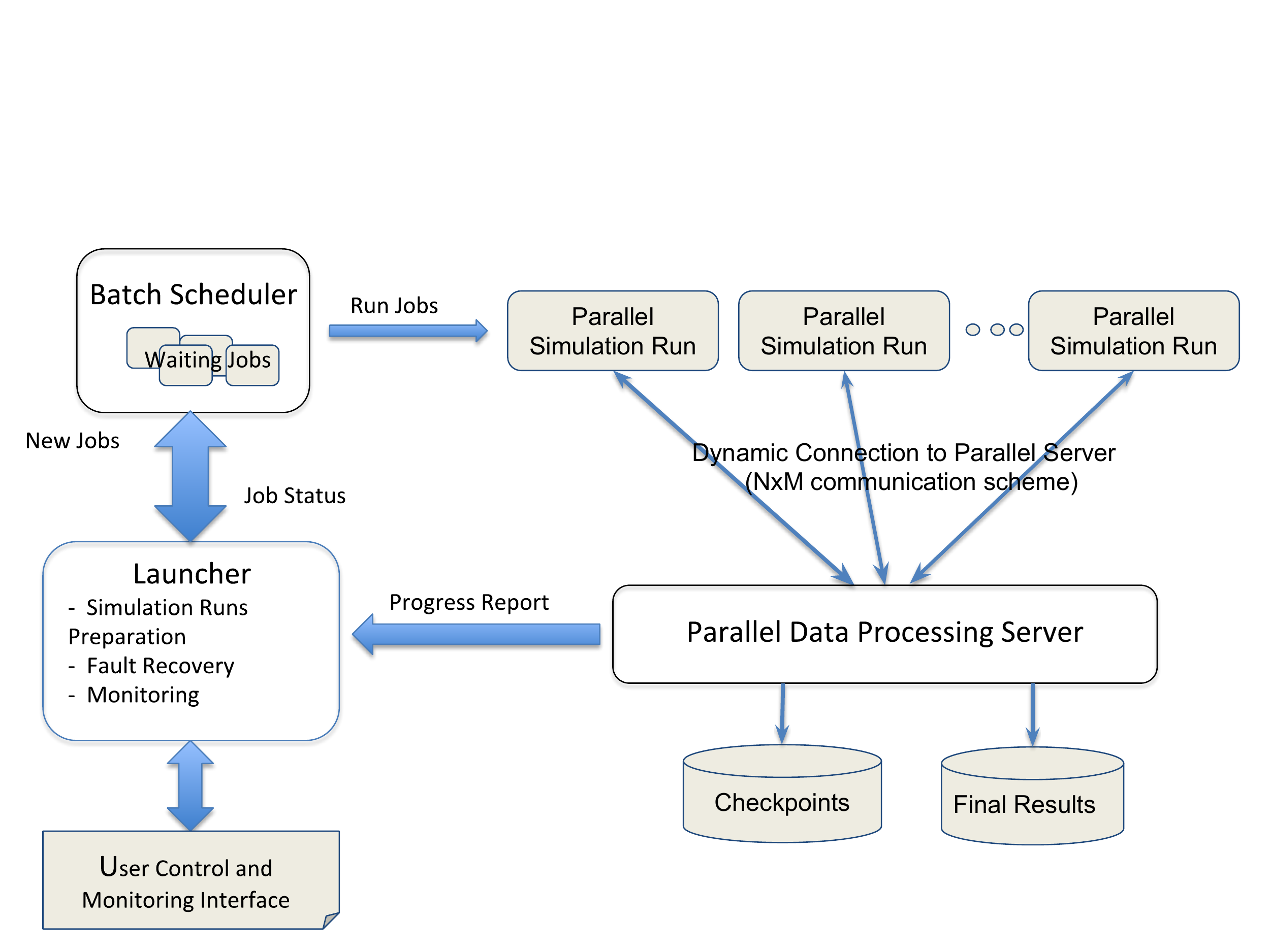}

\caption{ \melissa three tier architecture.  The launcher oversees the
execution in tight link with the batch scheduler. 
The job scheduler regulates the number of simulation jobs to run
according to the machine availability,  leading to an elastic resource
usage.   The parallel server, started first,  process incoming
data   as soon as  received from the  connected simulations.  
A fault tolerance mechanism automatically restarts failing 
simulation runs  or a failing parallel server.}
\label{fig:launcher}
\end{figure}

\melissa (Modular External Library for
  In Situ Statistical Analysis) is an open source
  framework~\footnote{\url{https://melissa-sa.github.io}}
  that relies on a three tier architecture  (Fig.~\ref{fig:launcher}). 
    The {\it \melissa clients}  are the
parallel simulations, providing their outputs to the server. The  {\it
  \melissa Server} aggregates the  simulation
results and updates iterative statistics as soon
as a new result is available.    {\it \melissa
 Launcher}  interacts with the supercomputer batch scheduler and 
\melissa Server,  for creating, launching, and supervising
the server and clients.

\subsubsection{\melissa Server}\label{sec:server}

Melissa Server is parallel and runs on several nodes.
The number of nodes required for the server is mainly defined by its memory needs.
The amount of memory needed  for each computed statistic field is in  the order of the size of the output field
of one simulation (number of time-steps $\times$  the  number of cells or points in the mesh).
The simulation domain is evenly  partitioned in space among the different server processes at starting time.
Melissa uses its own static space partitioning of the data. This
partitioning is different than the simulation partitioning, requiring
a data redistribution between each client and the server.
Updating the statistics is a local operation that requires neither communication nor synchronization  between the server processes.
The one-pass statistic algorithms allow the  server processes to  update their local statistics each time they receive a new data message, coming from any simulation, in any order.
%

\subsubsection{Dynamic Connection to \melissa Server}

When   a simulation starts, it  dynamically
connects to  \melissa Server. Each  simulation
process opens individual communication channels to
each server process that needs data   according to  the  data redistribution pattern.
Every time new results  are available, simulation processes send  the results toward \melissa Server. 

\melissa was designed to keep intrusion into the
simulation code minimal.    \melissa provides
three functions to integrate in the simulation code
through a dynamic library. The first
function (Initialize) allocates internal
structures and connects the simulation to the
server. At each time-step, the second function (Send) sends the simulation data to its corresponding \melissa Server processes.
The third function (Finalize) disconnects the simulation and releases the allocated structures. 

Melissa components are connected by \zmq communication sockets. 
  \zmq is a multi-threaded library for the efficient asynchronous transfer of
messages  between  a  client  and  a  server~\cite{zeromq-book-2013}.
\zmq bufferizes messages in a background thread both on the client and server side.
This allows to regulate data transfers between clients and server without blocking the executions.
Communications only become blocking when both buffers are full.

%


\subsubsection{\melissa Launcher}

\melissa Launcher  takes care of   generating   the
parameter sets, requesting the batch scheduler to start the server and
the clients, and track the various running job  progress.
It first submits  to the batch scheduler a job for the   \melissa Server.
Then, the launcher retrieves the server node addresses (the server is
parallelized on several nodes)  and submits  the simulation jobs.  
Next, each simulation is submitted to the batch scheduler in an independent
job, making Melissa very elastic.  Simulations can
be submitted all at once or at a more regulated pace depending on the
cluster policy for job submissions.

The launcher is the main actor of Melissa  fault tolerance
mechanism. It regularly checks for the jobs status, receives a heartbeat from the server,
and is able to resubmit the server or the simulation jobs if needed.

For each   use case,  the user  needs to  provide a
script for  the \melissa Launcher to  generate the
parameter  sets  and   to  launch  the  simulations.

\subsubsection{Fault Tolerance}

\melissa asynchronous client/server architecture
leverages the iterative statistics  computations 
to  support  a simple yet robust  fault tolerance
mechanism.  Melissa supports detection and recovery from
failures (including straggler issues) of \melissa Server
and simulations, through heartbeats and server check-pointing. 
Melissa Launcher  communicates  with the server and the batch scheduler
to detect simulation or server faults. As every simulation runs in a
separate job, the failure of one simulation does not impact the
ongoing study:  Melissa launcher simply restarts it and the server
discards already processed  messages. 

Melissa Server  regularly checkpoints. On failure,  Melissa Launcher 
kills  the running simulation jobs, restarts the server from the last
checkpoint and the associated  missing simulations. 

Errors on Melissa Launcher are fail-safe:  the running simulations
proceed to completion   with the server aggregating the incoming data.  
After a given time without any new message,  the server checkpoints
and stops. 

The server check-pointing enables to  manually restart
any   study for adding extra simulations. This is convenient 
for instance if  the system killed the running experiment because it
reached the wall-time limit, or simply because the user estimates that
more simulation runs are required for improving the quality of the
statistics.

%
%
%

\section{Experiments}\label{sec:experiments}

This section presents the large scale experiment illustrating the computation of  one-pass quantiles with Melissa.

%

\subsection{Fluid Simulation with \codesaturne}

%

The fluid numerical simulation is performed with \codesaturne\cite{saturne-ijfv}, an open-source computational fluid dynamics tool designed to solve the Navier-Stokes equations, with a focus on incompressible or dilatable flows  and advanced turbulence modeling.
\codesaturne relies on a finite volume discretization and allows the use of various mesh types, using an unstructured polyhedral cell model, allowing hybrid and non-conforming meshes.
The parallelization is based on a classical domain partitioning using MPI, with an optional second (local) level using OpenMP~\cite{Fournier-Saturne-2011}.

\subsection{Use Case}\label{sec:usecase}

\begin{figure}[!ht]
  \centering
  \includegraphics[width=0.6\linewidth]{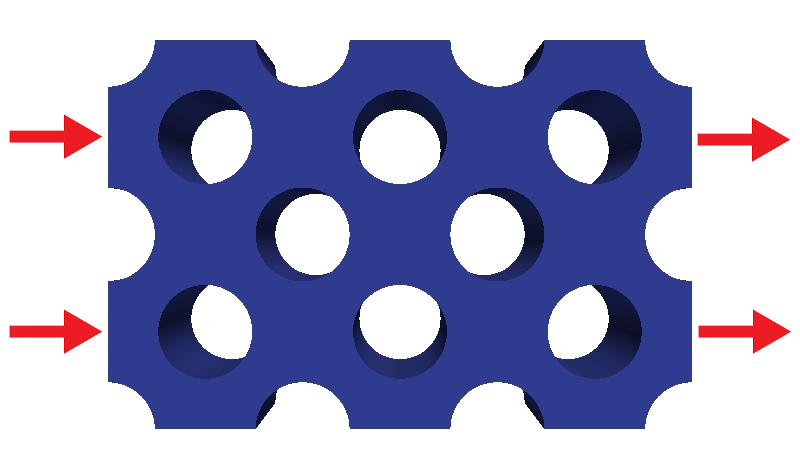}
  \vspace{-0.3cm}
  \caption{Use case: water flows from the left, between the tube bundle, and exits to the right.}
  \label{geo}
  \vspace{-0.3cm}
\end{figure}

Our implementation is validated on a fluid mechanics use case\hide{provided by EDF R\&D} simulating a water flow in a tube bundle (Fig. \ref{geo}).
The solved scalar field represents a dye concentration and could be replaced by temperature or concentration of chemical compounds in actual industrial studies.
The mesh is composed of $6 002 400$ hexahedra.
An ensemble  study is generated by simulating the injection of a tracer or dye along the inlet, with 2 independent injection surfaces, each defined by three varying parameters:
\begin{enumerate}
\item dye concentration on the upper inlet,
\item dye concentration on the lower inlet,
\item width of the injection on the upper inlet,
\item width of the injection on the lower inlet,
\item duration of the injection on the upper inlet,
\item duration of the injection on the lower inlet.
\end{enumerate}
All these parameters are uncertain and modelled as independent random variables.
Their probability distribution functions are uniform on $[0.002,0.1]$ for the two injection duration and uniform on $[0.1,0.9]$ for the two concentrations and the two injection widths.

To initialize our ensemble  study, we first ran a single $1000$
time-steps simulation, to obtain a steady flow. Assuming the resulting
flow is independent of the scalar (dye concentration) values, we then
use the final state of this simulation as the frozen velocity,
pressure, and turbulent variable fields, on which we perform our
experiment. This option allows solving only the convection-diffusion
equation associated to the scalar, so simulations run much faster
while generating the same amount of data. Each simulation consists of
100 time-steps on these frozen fields, with different parameter sets.

This study  ran a total  of $3 000$ simulations for computing all the ubiquitous 
percentiles, the 100-quantile, on the 6M  hexahedra and 100 time-steps.





\subsection{Performance}\label{sec:performances}

The experiment presented in this section run on a supercomputer called ``Eole'',  ranked 128th at  the top500.org of November 2016, when it was installed. At the time of submission of this article Eole ranks 460th, list   of November 2018. Eole is  composed of three kind of nodes:
\begin{enumerate}
\item 1164 standard nodes, each containing 2 Intel processors Xeon E5-2680v4 14C 2.4GHz and $128$ GB of memory. Each processor contains 14 cores thus a standard node contains 28 cores.
\item 162 big memory nodes with the same architecture than the standard nodes but with upgraded memory, between $256$ GB and $2$ TB.
\item 64 graphical nodes equipped with Nvidia K80 GPUs.
\end{enumerate}
All the nodes of Eole are connected by an Intel Omni-Path network.
In our experiment, each \codesaturne simulation runs on one node and it is parallelized on 28 cores.
On the server side, \melissa Server must have enough memory to keep all the updated statistics and to queue the inbound messages from the simulations, and must compute the statistics fast enough to consume the data faster than they arrive.
Otherwise, \melissa Server inbound message queue will end up full,  eventually  blocking  the simulations. For this study, \melissa Server runs on 8 nodes (224 cores).

We have run the same experiment on Occigen2, a supercomputer that ranks 77th at the top500.org of November 2018. This test was performed to evaluate robustness. Because our validation was positive and the produced results are the same in both supercomputers, we thus just present the experimentation performed on EOLE.

\begin{figure}[!ht]
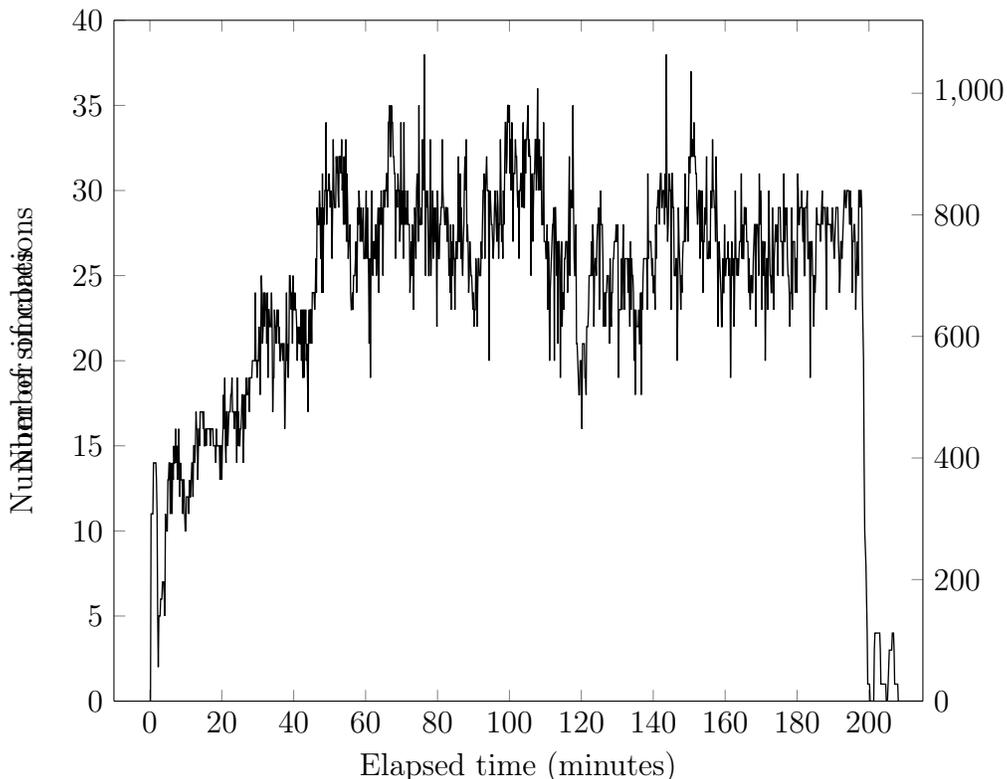

  \centering

\vspace{-0.2cm}
\vspace{-0.2cm}
\caption{Evolution of the number of simultaneously running simulations (equivalent to the number of cores) during the execution of our use case running a total of $3 000$ simulations.}
\label{fig:execution}
\end{figure}

Fig.~\ref{fig:execution} presents a plot showing the temporal evolution of the study. On the horizontal axis time evolves, in minutes, from left to right. Total execution time was 210 minutes (3.5 hours). The plot shows the evolution of the number of simultaneous fluid dynamics simulations performed over time. This is equivalent to the number of cores because each simulation runs, for this experiment, of 28 cores. We fixed a limit of 40 simultaneous simulations (1120 cores) but we see that this limit is seldom achieved because the scheduler of Eole allocates less resources for our study. This is a common situation faced by most users of supercomputers, which being shared resources are not available at will. We see how the elastic nature of Melissa allows the study to adapt to the available resources, sometimes using less resources, sometimes the full capacity (40 simultaneous simulations in this case) and even not running any simulation for some period of time, as we can see around time 200 minutes.

During this study, \melissa Server processed 11 TB of data coming
from the simulations. In a classical study, all these data would be
written to the file-system, and read back to compute the quantiles. This would not have been possible on Eole or Occigen simply due to the storage capacity (quota) being limited per user.

\section{Ubiquitous Quantile  Visualization}\label{sec:visu}

\begin{figure}[!ht]
\begin{subfigure}[t]{0.49\linewidth}
  \centering
  \captionsetup{width=.8\linewidth}
  \includegraphics[width=1\linewidth]{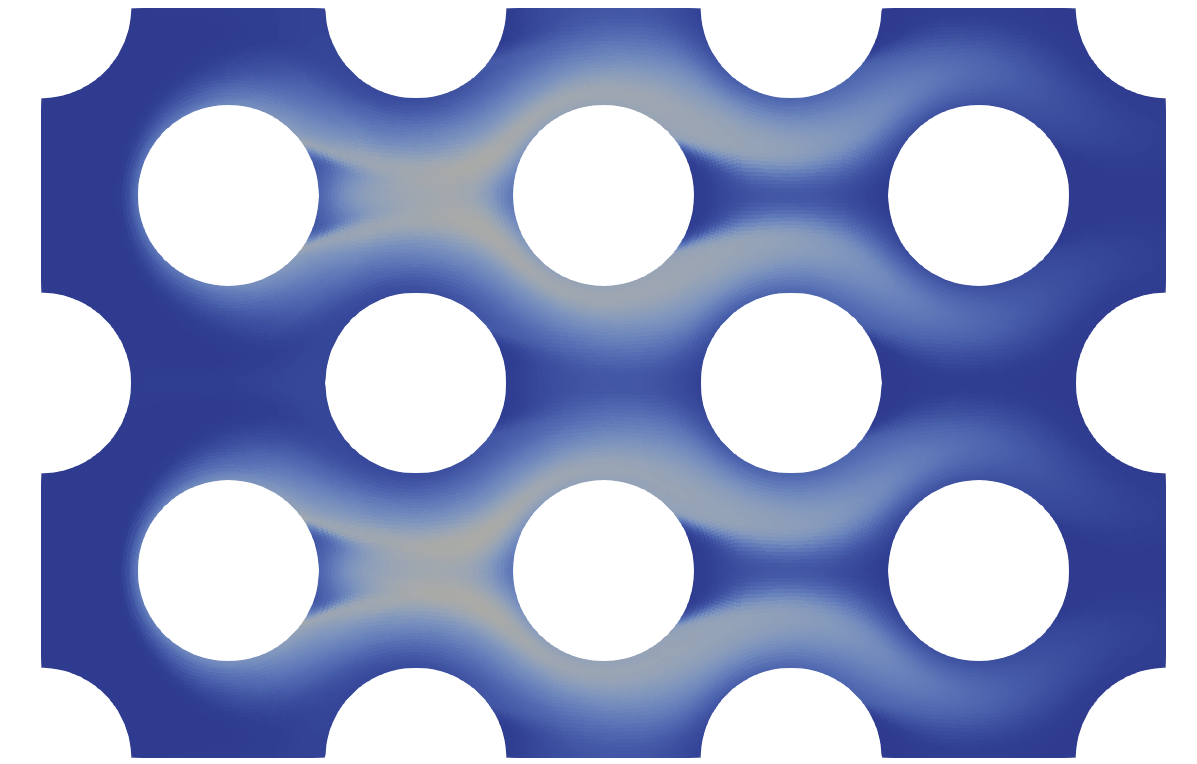}
  \caption{$75^{th}$ percentiles}
  \label{subfig:visutop1}
\end{subfigure}
\begin{subfigure}[t]{0.49\linewidth}
  \captionsetup{width=.8\linewidth}
  \includegraphics[width=1\linewidth]{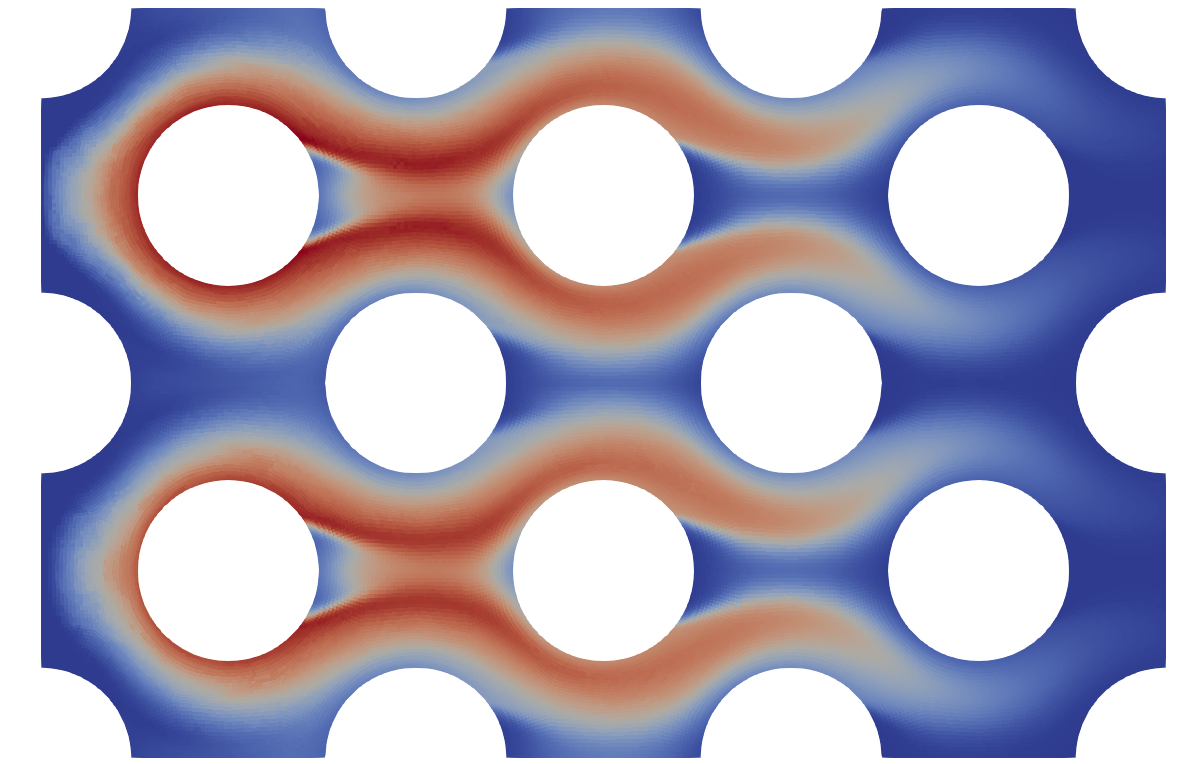}
  \caption{$95^{th}$ percentiles}
  \label{subfig:visutop2}
\end{subfigure}

\begin{subfigure}[t]{0.49\linewidth}
  \centering
  \includegraphics[width=1\linewidth]{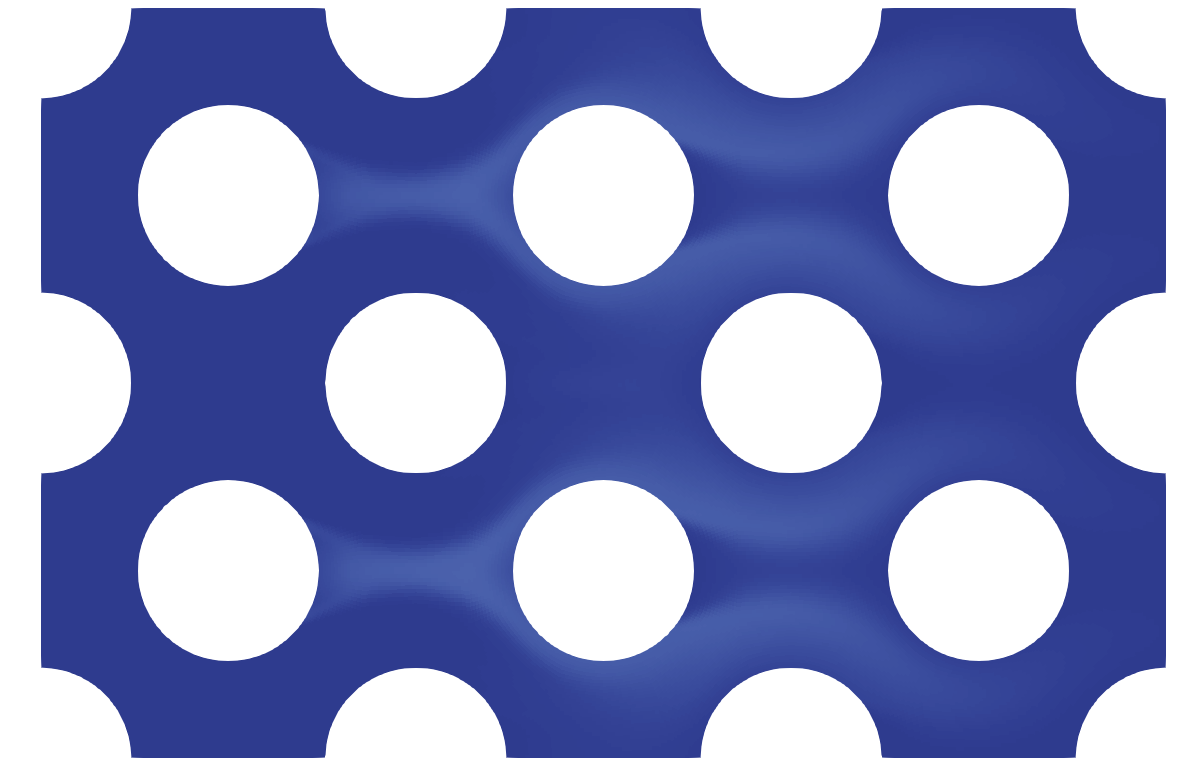}
  \caption{$25^{th}$ percentiles}
  \label{subfig:visumiddle1}
\end{subfigure}
\begin{subfigure}[t]{0.49\linewidth}
  \includegraphics[width=1\linewidth]{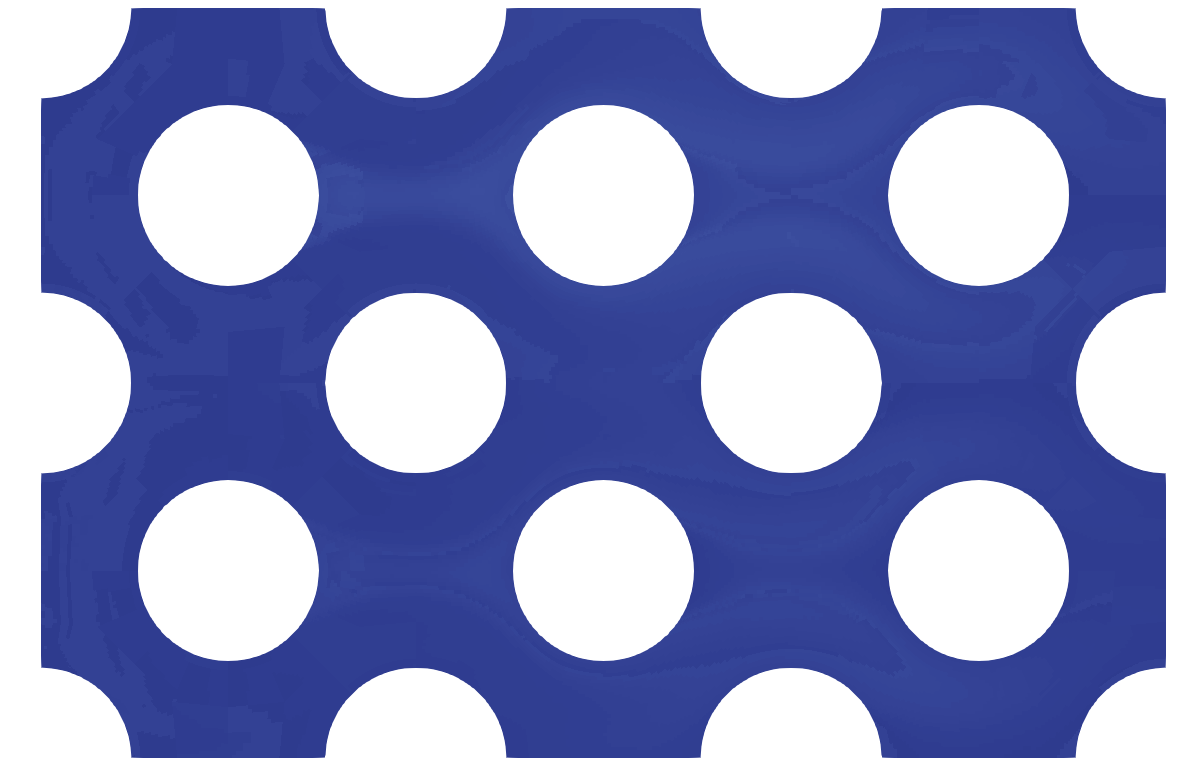}
  \caption{$5^{th}$ percentiles}
  \label{subfig:visumiddle2}
\end{subfigure}

\begin{subfigure}[t]{0.49\linewidth}
  \centering
  \includegraphics[width=1\linewidth]{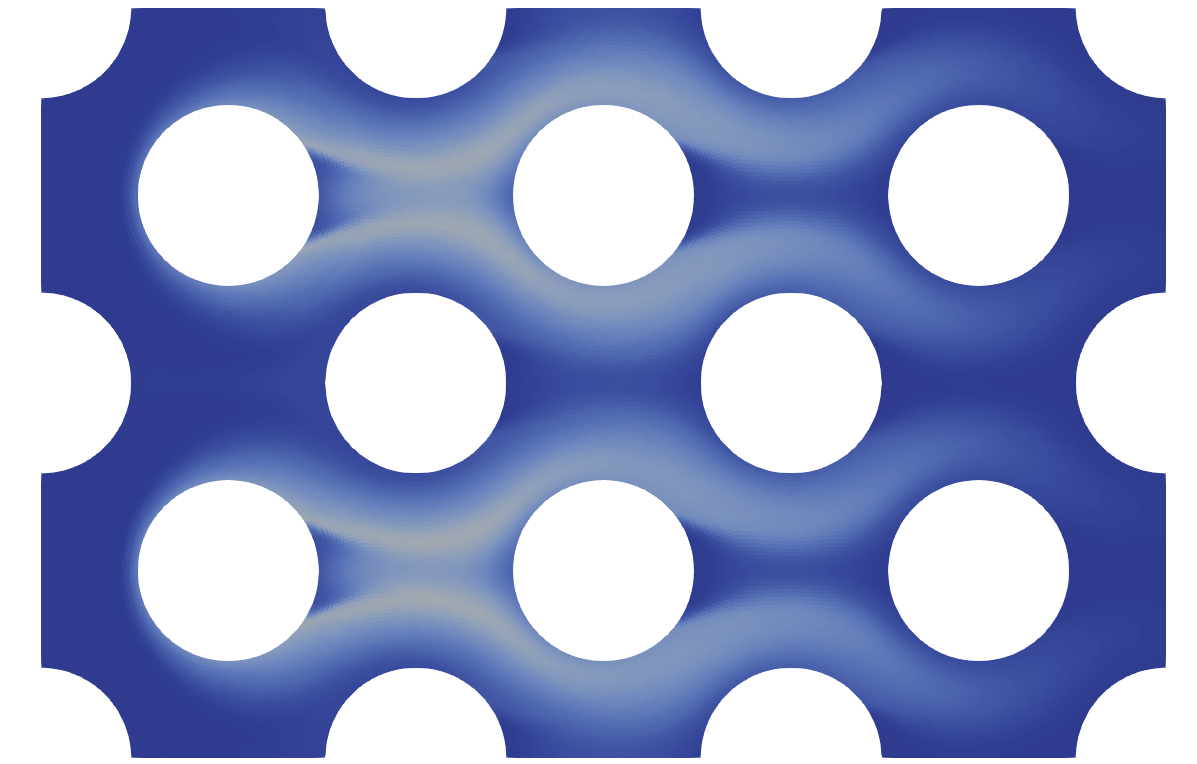}
  \caption{50\% inter-percentile range}
  \label{subfig:visubottom1}
\end{subfigure}
\begin{subfigure}[t]{0.49\linewidth}
  \includegraphics[width=1\linewidth]{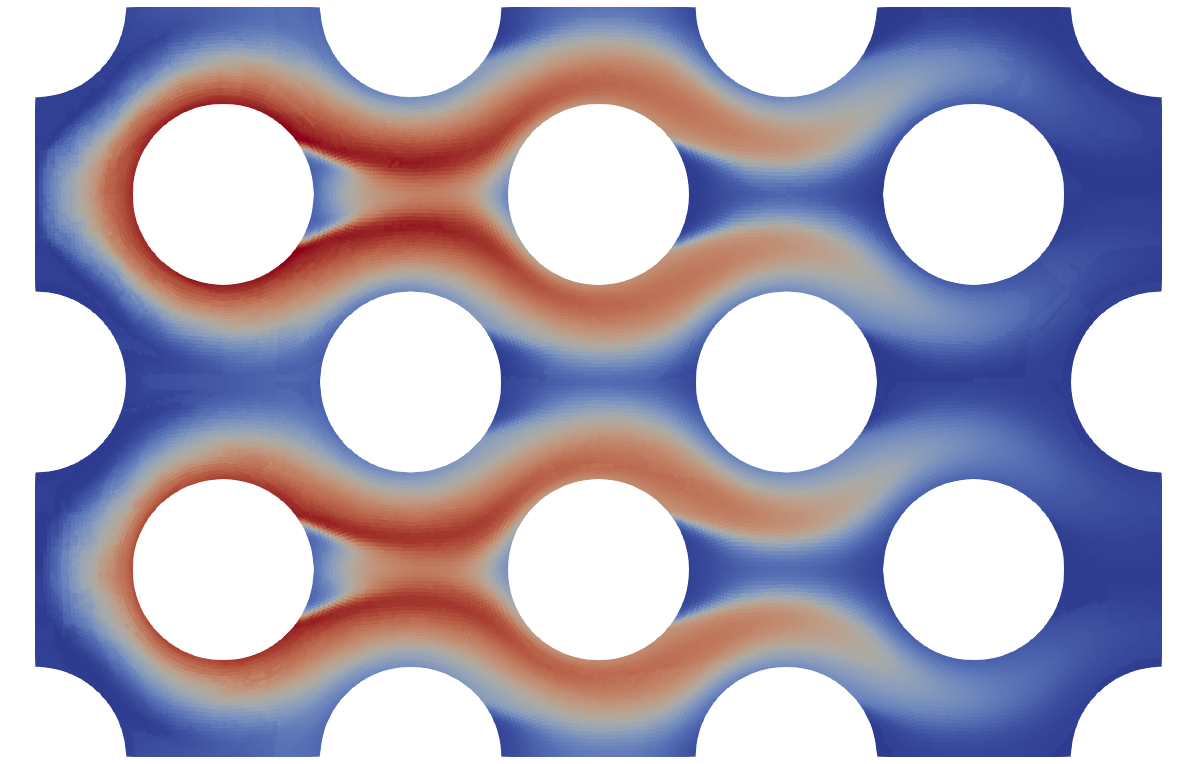}
  \caption{90\% inter-percentile range}
  \label{subfig:visubottom2}
\end{subfigure}

\vspace{0.1cm}
  \caption{Percentiles and inter-percentile range maps on a slice of the mesh at $80^{th}$ time-step.  The four top panels correspond to the percentiles while the two bottom panels correspond to inter-percentile ranges.
   All maps share the same scale.} 
\vspace{-0.3cm}
  \label{fig:visuquantiles}
\end {figure}

This section presents the quantiles computed during the experiments. 
Figure~\ref{fig:visuquantiles} presents six spatial maps extracted from the ubiquitous quantiles. By use of the Open-Source visualization tool ParaView, we have chosen a time-step and performed a slice on a mid-plane aligned with the direction of the fluid. The chosen time-step belongs to the last temporal part of the simulation ($80^{th}$ time-step over $100$). 
This section focuses on the interpretation of the computed percentiles. However, the system can compute any other kind of quantile.

On the four top panels of Figure~\ref{fig:visuquantiles}, Fig.~\ref{subfig:visutop1}, \ref{subfig:visutop2}, \ref{subfig:visumiddle1} and \ref{subfig:visumiddle2}, we present the $75^{th}$, $95^{th}$, $25^{th}$, $5^{th}$ percentiles, respectively. 
On the two bottom panels (Fig.~\ref{subfig:visubottom1} and Fig.~\ref{subfig:visubottom2}) the inter-percentile ranges containing 50\% and 90\% of the samples are shown. 
Inter-percentile ranges are easily computed from percentiles by substraction: the 50\% inter-percentile range corresponds to the $75^{th}$ percentile minus the $25^{th}$ percentile; the 90\% inter-percentile range corresponds to the $95^{th}$ percentile minus the $5^{th}$ percentile. 
In Figure~\ref{fig:visuquantiles} each column shows an inter-percentile map on the bottom and the percentile maps that served for its calculation above it. Looking at these maps an analyst can deduce several things:
\begin{enumerate}
\item Extreme percentile maps such as $95^{th}$,
  Fig.~\ref{subfig:visutop2}, give an idea of the distribution of the
  upper bounds of all simulations. In our use case, we can assess
  which spatial areas contain low quantities of dye. Areas
  colored in blue necessarily contain low dye concentrations for any
  simulation in the ensemble study. Extreme low percentile maps, such as $5^{th}$ has also a direct interpretation in the opposite sense. 
\item Inter-percentile range maps such as Fig.~\ref{subfig:visubottom1} or ~\ref{subfig:visubottom2} are maps that show the spatial variability of statistical dispersion. Indeed, scalar inter-percentile ranges are non-parametric measures of statistical dispersion, which means that no a priori knowledge about the distribution of the data is needed. This characteristic makes these ranges both general and robust. Visualizing a map of such a measure of dispersion allows to understand how the data distribution spatially concentrates. In our use case, the low percentile maps used to calculate the inter-percentile maps are mainly close to zero for all cells of the mesh, which makes these maps resemble the higher percentiles maps. However, this is in general not true.
\end{enumerate}

\begin {figure}[!ht]
\begin{subfigure}[t]{0.9\linewidth}
  \centering
  \captionsetup{width=.8\linewidth}
  \includegraphics[width=0.6\linewidth]{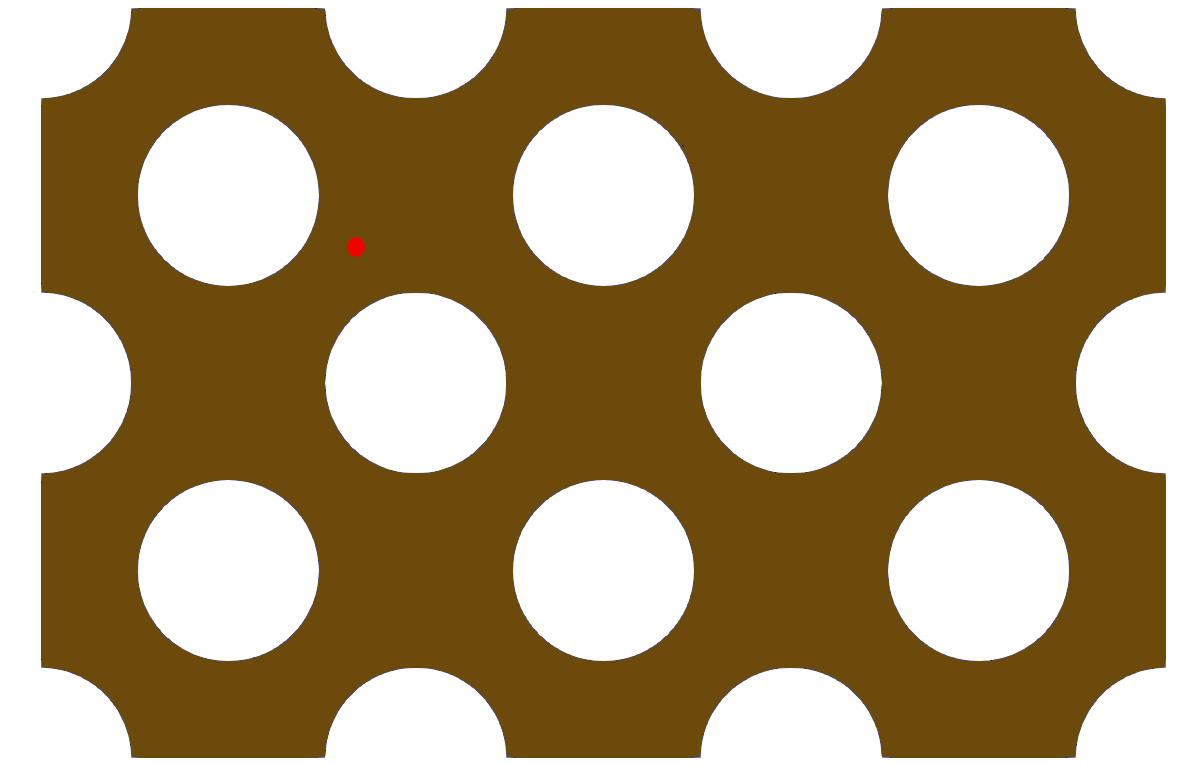}
  \caption{Probe position}
  \label{subfig:timeup}
\end{subfigure}

\begin{subfigure}[t]{0.98\linewidth}
  \centering
  \includegraphics[width=0.8\linewidth]{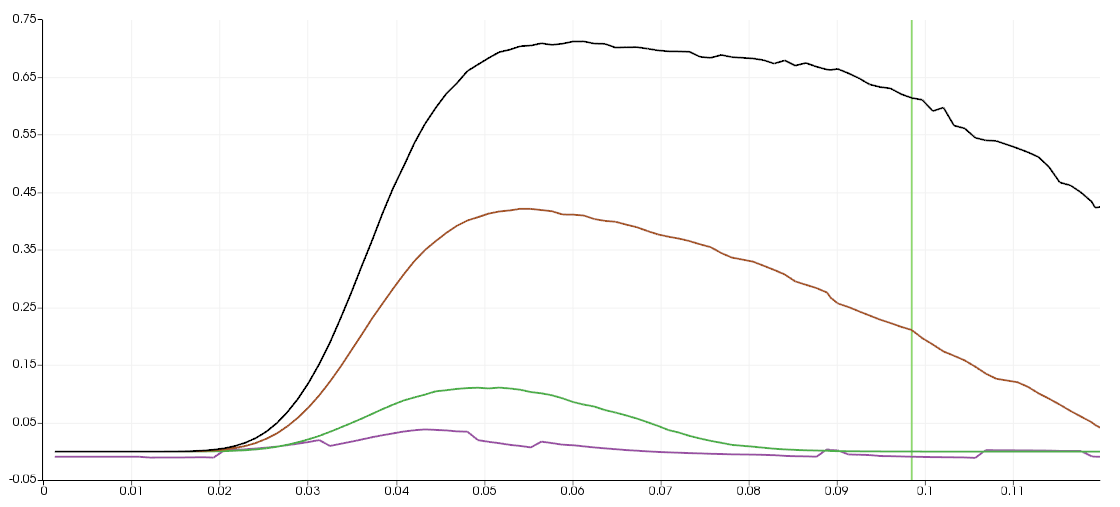}
  \caption{Evolution of the percentiles for the dye concentration (vertical axis) over time (horizontal axis)}
  \label{subfig:timebottom}
\end{subfigure}

\vspace{0.1cm}
  \caption{A probe in a cell of the mesh allows an extraction of the temporal evolution of percentiles at a specific spatial location.} 
\vspace{-0.3cm}
  \label{fig:timeevolution}
\end {figure}

The maps shown in Figure~\ref{fig:visuquantiles} are static and 2D but
we recall that we calculate ubiquitous percentiles, thus 3D and time
dependent data is available. 
The Figure~\ref{fig:timeevolution} shows the temporal
evolution of a probe positioned in the mesh using ParaView. At a
specific location, a temporal evolution of all computed quantiles can
be performed. In Fig.~\ref{subfig:timebottom} this evolution is plotted
for the  $95^{th}$, $75^{th}$, $25^{th}$ and $5^{th}$ percentiles. The
vertical line indicates the position of the current time step
($80^{th}$ time step). This figure clearly shows how the output
variability of the ensemble study depends on time. Indeed, all simulations contain no dye for the first 15 time-steps, which is the time the dye takes to propagate from the top injector to the spatial location of the probe. After this point, we observe a moment where the variability of the dye concentration is the highest before  a general decrease. 

Fig.~\ref{subfig:timebottom} can be seen as the evolution of a Tukey boxplot~\cite{tukey1977exploratory} over time. In fact, the  $25^{th}$ and $75^{th}$ percentiles correspond to the $1^{st}$ and $3^{rd}$ quartiles thus delimit the central box of the plot, while the $5^{th}$ and $95^{th}$ quantiles can be a choice for the whiskers. Using this analogy, we can easily observe that the dispersion of the dye concentration on the whole ensemble moves over time. Furthermore, the distribution of this quantity is not symmetrical and its asymmetry is evolving over time. 

\begin{figure}[!ht]
  \centering
  \includegraphics[width=0.8\linewidth]{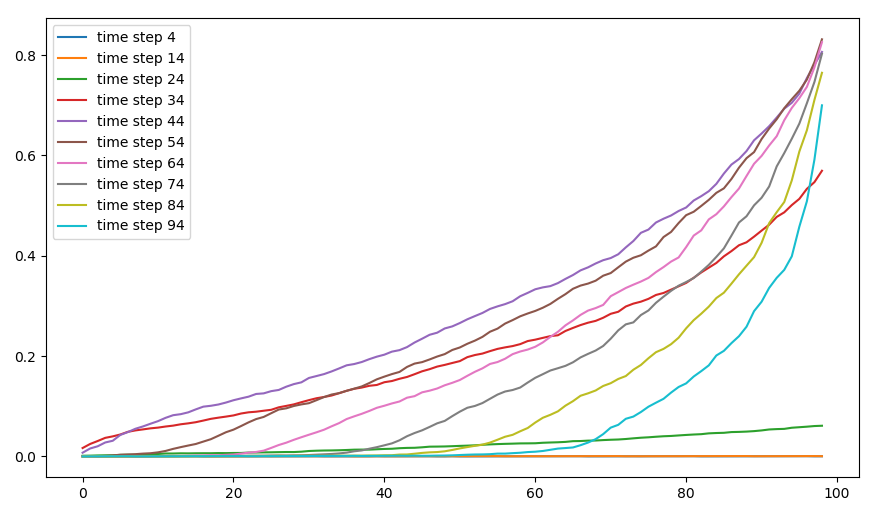}
  \vspace{-0.3cm}
  \caption{Percentile functions of the dye concentration at different time steps of the simulation. Vertical axis represents dye concentration. Horizontal axis represents percentiles. Each curve corresponds to different time steps of the simulation. All curves have been extracted form the probe position shown in Fig.~\ref{subfig:timeup}.} 
  \vspace{-0.3cm}
  \label{fig:quantiles}
\end{figure}

Finally, Fig.~\ref{fig:quantiles} shows a different representation of the evolution of the dye concentration at a fixed probe (the same than in Fig.~\ref{fig:timeevolution}).
At different regularly sampled time steps, the quantile functions of the concentration values are plotted (as a function of the order of the quantiles, between $0\%$ and $100\%$).
One can first observe zero-valued quantile functions for the first time steps (time steps 4 and 14).
Indeed, at the probe, the dye concentration is zero during the first times of the injection. 
Then, from time step 24 to time step 44, all the values of the quantile functions regularly increase.
It means that the dye concentration values homogeneously increase from $0$, reaching a maximal value close to $0.82$ for the $100\%$-order quantile. 
At the end of the simulation time, from time step 54 to time step 94, the quantile functions are regularly displaced to the right.
The values close to zero disappear and the concentration of strong values becomes more and more important.
As a conclusion, thanks to the the quantile functions, this graph  allows to finely and quantitatively analyze the temporal evolution of this dye concentration phenomena.
We remark that, because we have calculated the ubiquitous percentiles, it is possible to obtain Fig.~\ref{fig:timeevolution} and Fig.~\ref{fig:quantiles} for any location on the simulation domain.

\section{Conclusion}\label{sec:conclusion}

Exascale machines will shortly become a reality. But so far only very
few applications are able to  take benefit of the level
of parallelization these machines provide.  Ensemble runs and uncertainty quantification approaches may require to execute from thousands to millions of the same simulation,
 making it an extremely compute-intensive process that will fully
 benefit from Exascale machines. However, the large amount of data
 generated is a strong  I/O bottleneck if these intermediate data are
 saved to disk. An alternative consists in processing these data in
 transit as proposed by the Melissa framework.   But this approach
requires one-pass data processing algorithms to limit the amount of
memory needed.  

This paper  proposes the computation of quantiles by use of a parallel one-pass strategy based on a new robust version of the stochastic quantile algorithm of Robbins-Monro~\cite{robmon51}. 
The algorithm is experimented at large scale with the Melissa in transit, elastic and fault-tolerant  processing framework. 
On a fluid dynamics application case, $3 000$ simulations on a $6M$ cell mesh allow to compute all spatio-temporal percentiles at full resolution, saving $11$TB of intermediate storage thanks to Melissa. 
Because the I/O bottleneck has been bypassed, ubiquitous spatio-temporal maps of percentiles and inter-percentile based intervals have been successfully visualized, revealing  their interest for the interpretation by the users.
In addition to quantiles, Melissa  currently supports average, standard deviation, skewness,  kurtosis, minimum, maximum, threshold exceedance and Sobol' indices.

Several ways of improvement are identified on the iterative quantile estimation algorithm that we use. 
First, the constant $C$ and the $\gamma$-profile has been defined in an heuristic manner and a theoretical consolidation of such a choice seems necessary.
Different versions of the Robbins-Monro algorithm also exist (for
example the averaged one) and should be tested. Second, the user of
the method should avoid to fix a priori a total number of runs. A more
interesting approach would iteratively give  an estimate of the
quantile, and then stop when a sufficient precision (defined by the
user) is achieved. This precision is not theoretically known in a non-asymptotic context but approximate results are maybe accessible and have to be tested. 

Last, the quantiles of the spatio-temporal outputs have been
computed cell per cell and time-step per time-step.
The interpretation of this ubiquitous quantiles 
(for instance in the form of static spatial maps, temporal probes or videos) is
 much richer than the traditional predefined probe-based or sub-sampled approaches.
However, the functional space where the spatio-temporal field lives has not be estimated.
Dealing with this space (as in \cite{carcen13}), the ubiquitous quantile estimates would conserve the geometrical and temporal
structure of the ensemble run  study.

\section{Acknowledgement}

This work was partly funded in the Programme
d'Investissements d'Avenir, grant PIA-FSN2-Calcul intensif et simulation num\'erique-2-AVIDO.
This  work  was  granted  access  to  the  HPC  resources  of CINES under  the allocation 2018-A0040610366 made by GENCI.



\section{References}








\end{document}